\documentclass[a4paper,12pt]{article}

\usepackage{amsmath,amssymb,amsthm,amscd}

\newcommand{\Fb}{\ensuremath{\mathfrak{b}}}
\newcommand{\Fg}{\ensuremath{\mathfrak{g}}}
\newcommand{\Fh}{\ensuremath{\mathfrak{h}}}
\newcommand{\Fn}{\ensuremath{\mathfrak{n}}}

\newcommand{\BQ}{\ensuremath{\mathbb{Q}}}

\newcommand{\ch}{\mathop{\rm ch}\nolimits}
\newcommand{\tch}{\mathop{\rm ch}^{\omega}\nolimits}
\newcommand{\tr}{\mathop{\rm tr}\nolimits}

\newcommand{\pos}{P_{\omega}^{\ast}}

{\begin{list}{}{\setlength{\topsep}{2mm}\setlength{\itemsep}{-1mm}}}%
{\end{list}}%

\newcommand{\ha}[1]{\ensuremath{\widehat{#1}}}
\newcommand{\ti}[1]{\ensuremath{\widetilde{#1}}}

\makeatletter
\renewcommand\section{\@startsection{section}{1}{0pt}
{2.5ex plus 1ex minus .2ex}{1.0ex plus .2ex}{\large\sc}}
\renewcommand\subsection{\@startsection{subsection}{1}{0pt}
{2.5ex plus 1ex minus .2ex}{-1em}{\bf}}

\theoremstyle{plain}

\newtheorem{introthm}{Theorem}

\theoremstyle{definition}

\theoremstyle{remark}

\begin{document}
\setlength{\baselineskip}{18pt}
\setcounter{section}{-1}


\title{{\sc A Twining Character Formula for Demazure Modules}}
\author{By \\[2mm] Daisuke Sagaki \\ 
\small Graduate School of Mathematics, 
University of Tsukuba, \\
\small Tsukuba, Ibaraki 305-8571, Japan \ 
(e-mail: sagaki@math.tsukuba.ac.jp)}
\date{}
\maketitle


In \cite{FRS} and \cite{FSS}, they introduced new character-like
quantities corresponding to a graph automorphism of a Dynkin diagram, 
called twining characters, for certain Verma modules and integrable
highest weight modules over a symmetrizable Kac-Moody algebra, 
and gave twining character formulas for them. 
Recently, the notion of twining characters has naturally been
extended to various modules, and formulas for them 
has been given (\cite{KN}, \cite{KK}, \cite{N1}--\cite{N4}). 
In this short note, we introduce a twining character formula 
for Demazure modules over a symmetrizable 
Kac-Moody algebra, which is obtained in \cite{S}.

\vspace{3mm}

Let $\Fg=\Fg(A)=\Fn_{-} \oplus \Fh \oplus \Fn_{+}$ 
be a symmetrizable Kac-Moody algebra over $\BQ$ associated 
to a GCM $A=(a_{ij})_{i,j \in I}$ of finite size, 
where $\Fh$ is the Cartan subalgebra, $\Fn_{+}$ the sum of positive root
spaces, and $\Fn_{-}$ the sum of negative root 
spaces, and let $\omega:I \, \rightarrow \, I\,$ 
be a (Dynkin) diagram automorphism, that is, a bijection 
$\omega:I \, \rightarrow \, I$ satisfying 
$a_{\omega(i),\,\omega(j)}=a_{ij}$ for all $i,j \in I$.
It is known that a diagram automorphism induces 
a Lie algebra automorphism $\omega \in \text{Aut}(\Fg)$ 
that preserves the triangular decomposition of $\Fg$ 
(see \cite{FSS} and \cite{FRS}). 
Then we define a linear automorphism 
$\omega^{\ast} \in \text{GL}(\Fh^{\ast})$ by $(\omega^{\ast}(\lambda))(h):=
\lambda(\omega(h))$ for $\lambda \in \Fh^{\ast}$, $h \in \Fh$. 
We set $(\Fh^{\ast})^{0}:=\{\lambda \in \Fh^{\ast} \, | \, 
\omega^{\ast}(\lambda)=\lambda\}$, and 
call its elements symmetric weights. 
We also set $\ti{W}:=\{w \in W \, | \, 
w\,\omega^{\ast}=\omega^{\ast}\,w\}$. 

Further we define a ``folded'' 
matrix $\ha{A}$ associated to $\omega$, which is 
again a symmetrizable GCM if $\omega$ satisfies a certain condition, 
called the linking condition (see \cite{FSS}). 
Then the Kac-Moody algebra $\ha{\Fg}=\Fg(\ha{A})$ 
associated to $\ha{A}$ is 
called the orbit Lie algebra. We denote by $\ha{\Fh}$ 
the Cartan subalgebra of $\ha{\Fg}$ and by $\ha{W}$ 
the Weyl group of $\ha{\Fg}$. 
Then there exist a linear isomorphism 
$\pos:\ha{\Fh}^{\ast} \, \rightarrow \, (\Fh^{\ast})^{0}$ and 
a group isomorphism $\Theta : \ha{W} \, \rightarrow \, \ti{W}$
such that $\Theta(\ha{w})=\pos \circ \ha{w} \circ (\pos)^{-1}$
for all $\ha{w} \in \ha{W}$.

Let $\lambda$ be a dominant integral weight. 
Denote by $L(\lambda)=\bigoplus_{\chi \in \Fh^{\ast}}
L(\lambda)_{\chi}$ the irreducible highest weight 
$\Fg$-module of highest weight $\lambda$. Then, for $w \in W$, 
we define the Demazure module $L_{w}(\lambda)$ of lowest weight 
$w(\lambda)$ in $L(\lambda)$ by $L_{w}(\lambda):=U(\Fb)u_{w(\lambda)}$, 
where $u_{w(\lambda)} \in L(\lambda)_{w(\lambda)} 
\setminus \{0\}$ and $U(\Fb)$ is the universal enveloping algebra of 
the Borel subalgebra $\Fb:=\Fh \oplus \Fn_{+}$ of $\Fg$. 
If $\lambda$ is symmetric, then we know from \cite{FSS} and \cite{FRS} 
(see also \cite{N1}) that there exists a (unique) linear automorphism 
$\tau_{\omega}:L(\lambda) \, \rightarrow \, L(\lambda)$ such that 
\begin{equation*}
\tau_{\omega}(xv)=\omega^{-1}(x)\tau_{\omega}(v) \quad 
\text{for all } x \in \Fg, \, v \in L(\lambda)
\end{equation*}
and $\tau_{\omega}(u_{\lambda})=u_{\lambda}$ with 
$u_{\lambda}$ a (nonzero) highest weight 
vector of $L(\lambda)$. Then it is easily seen 
that the Demazure module $L_{w}(\lambda)$ with $w \in \ti{W}$ 
is $\tau_{\omega}$-stable. Here we define 
the twining character $\tch (L_{w}(\lambda))$ 
of $L_{w}(\lambda)$ by:
\begin{equation*}
\tch (L_{w}(\lambda)):=\sum_{\chi \in (\Fh^{\ast})^{0}} 
\tr\bigl(\tau_{\omega}|_{L_{w}(\lambda)_{\chi}}\bigr)e(\chi).
\end{equation*}
Our formula is the following:
\begin{introthm}\ 
Let\, $\lambda$ be a symmetric dominant integral weight and\, 
$w \in \ti{W}$. Set\, $\ha{\lambda}:=(\pos)^{-1}(\lambda)$
and\, $\ha{w}:=\Theta^{-1}(w)$. Then we have
\begin{equation*}
\tch (L_{w}(\lambda))=\pos(\ch \ha{L}_{\ha{w}}(\ha{\lambda})),
\end{equation*}
where\, $\ha{L}_{\ha{w}}(\ha{\lambda})$ is the Demazure module 
of lowest weight\, $\ha{w}(\ha{\lambda})$ in the irreducible 
highest weight module\, $\ha{L}(\ha{\lambda})$ of highest weight 
$\ha{\lambda}$ over the orbit Lie algebra\, $\ha{\Fg}$.
\end{introthm}


\setlength{\baselineskip}{18pt}
\end{document}